\documentclass[11pt]{article}
\usepackage{amsmath,amssymb,theorem}
\usepackage{graphicx, enumerate}
\usepackage{subfigure}
\usepackage{psfrag}
\usepackage{placeins}
\usepackage{wrapfig}
\usepackage{booktabs}
\usepackage{chngcntr}
\counterwithin{figure}{section}
\numberwithin{figure}{section}
\counterwithin{table}{section}
\usepackage{hyperref}
\hypersetup{
  colorlinks,
  citecolor=black,
  filecolor=black,
  linkcolor=black,
  urlcolor=black
}

\newcommand{\mc}{\mathcal}

\newtheorem{thm}{Theorem}[section]
\newtheorem{conj}[thm]{Conjecture}

\newtheorem{claim}[thm]{Claim}
\newtheorem{lem}[thm]{Lemma}
\newtheorem{prop}[thm]{Proposition}
\theorembodyfont{\rmfamily}
\def\pf{\bigskip\noindent {\emph{Proof}.}~~}

\def\dfn#1{{\sl #1}}
\def\es{\emptyset}
\def\less{\setminus}

\def\qed{ \hfill $\blacksquare$}

\begin{document}
\title{Multicolor Gallai-Ramsey numbers of  $C_9$ and $C_{11}$}
\author{
Christian Bosse and Zi-Xia Song\thanks{Corresponding Author.    Email address:  Zixia.Song@ucf.edu.}\\
Department  of Mathematics\\
University of Central Florida\\
Orlando, FL 32816, USA\\
}
\maketitle
\begin{abstract}

A {\it Gallai coloring} is a coloring of the edges of a complete graph without rainbow triangles, and a {\it Gallai $k$-coloring} is a Gallai coloring that uses $k$ colors.   We  study Ramsey-type problems in Gallai colorings.  Given  an integer $k\ge1$ and a graph $H$, the  Gallai-Ramsey number $GR_k(H)$ is the least positive  integer $n$ such that every Gallai $k$-coloring of   the complete graph on $n$ vertices contains a monochromatic copy of  $H$.  It turns out that $GR_k(H)$ is more well-behaved than the classical Ramsey number $R_k(H)$. However,  finding exact values of $GR_k (H)$ is  far from trivial.  In this paper, we study Gallai-Ramsey numbers of odd cycles. We prove that for  $n\in\{4,5\}$  and all  $k\ge1$, $GR_k(C_{2n+1})= n\cdot 2^k+1$.  This  new result provides partial evidence for the first two open cases of the Triple Odd Cycle Conjecture of  Bondy and Erd\H{o}s from 1973. 
Our technique relies heavily on the structural result of Gallai on Gallai colorings of complete graphs. We believe the method we developed can be used to  determine the exact values of $GR_k(C_{2n+1})$ for all $n\ge6$.

\end{abstract}

\section{Introduction}
 \baselineskip 16pt

All graphs in this paper are finite and simple; that is, they have no loops or parallel edges. Given a graph $G$ and a set $A\subseteq V(G)$,  we use   $|G|$    to denote  the  number of vertices    of $G$, and  $G[A]$ to denote the  subgraph of $G$ obtained from $G$ by deleting all vertices in $V(G)\less A$.  A graph $H$ is an \dfn{induced subgraph} of $G$ if $H=G[A]$ for some $A\subseteq V(G)$.  We use $K_n$ and  $C_n$   to denote the      complete graph and cycle  on $n$ vertices, respectively.  
For any positive integer $k$, we write  $[k]$ for the set $\{1,2, \ldots, k\}$. We use the convention   ``$A:=$'' to mean that $A$ is defined to be the right-hand side of the relation.

 Given an integer $k \ge 1$ and a graph  $H $, the classical Ramsey number $R(H)$   is  the least    integer $n$ such that every $k$-coloring of  the edges of  $K_n$  contains  a monochromatic copy of  $H$.  Ramsey numbers are notoriously difficult to compute in general. In this paper, we  study Ramsey numbers of graphs in Gallai colorings, where a \dfn{Gallai coloring} is a coloring of the edges of a complete graph without rainbow triangles (that is, a triangle with all its edges colored differently). Gallai colorings naturally arise in several areas including: information theory~\cite{KG}; the study of partially ordered sets, as in Gallai's original paper~\cite{Gallai} (his result   was restated in \cite{Gy} in the terminology of graphs); and the study of perfect graphs~\cite{CEL}. There are now a variety of papers  which consider Ramsey-type problems in Gallai colorings (see, e.g., \cite{chgr, c5c6,GS, exponential, Hall, DylanSong,  C6C8, C10}).   These works mainly focus on finding various monochromatic subgraphs in such colorings. More information on this topic  can be found in~\cite{FGP, FMO}.  
 
A \dfn{Gallai $k$-coloring} is a Gallai coloring that uses $k$ colors. 
 Given an integer $k \ge 1$ and a graph  $H$, the   \dfn{Gallai-Ramsey number} $GR_k(H)$ is the least integer $n$ such that every Gallai $k$-coloring of $K_n$   contains a monochromatic copy of $H$.    Clearly, $GR_k(H) \leq R_k(H)$ for all $k\ge1$ and $GR_2(H) = R_2(H)$.    In 2010, 
Gy\'{a}rf\'{a}s,   S\'{a}rk\"{o}zy,  Seb\H{o} and   Selkow~\cite{exponential} proved   the general behavior of $GR_k(H)$.

\begin{thm} [\cite{exponential}]
Let $H$ be a fixed graph  with no isolated vertices 
 and let $k\ge1$ be an integer. Then
$GR_k (H) $ is exponential in $k$ if  $H$ is not bipartite,    linear in $k$ if $H$ is bipartite but  not a star, and constant (does not depend on $k$) when $H$ is a star.			
\end{thm}

It turns out that for some graphs $H$ (e.g., when $H=C_3$),  $GR_k(H)$ behaves nicely, while the order of magnitude  of $R_k(H)$ seems hopelessly difficult to determine.  It is worth noting that  finding exact values of $GR_k (H)$ is  far from trivial, even when $|H|$ is small.
We will utilize the following important structural result of Gallai~\cite{Gallai} on Gallai colorings of complete graphs.

\begin{thm}[\cite{Gallai}]\label{Gallai}
	For any Gallai-coloring $c$ of a complete graph $G$, $V(G)$ can be partitioned into nonempty sets  $V_1, V_2, \dots, V_p$ with $p>1$ so that    at most two colors are used on the edges in $E(G)\less (E(V_1)\cup \cdots\cup  E(V_p))$ and only one color is used on the edges between any fixed pair $(V_i, V_j)$ under $c$, where $E(V_i)$ denotes the set of edges in $G[V_i]$ for all $i\in [p]$. 
\end{thm}

The partition given in Theorem~\ref{Gallai} is  a \dfn{Gallai-partition} of  the complete graph $G$ under  $c$.  Given a Gallai-partition $V_1, V_2, \dots, V_p$ of the complete graph $G$ under $c$, let $v_i\in V_i$ for all $i\in[p]$ and let $\mathcal{R}:=G[\{v_1, v_2, \dots, v_p\}]$. Then $\mathcal{R}$ is  the \dfn{reduced graph} of $G$ corresponding to the given Gallai-partition under $c$. Clearly,  $\mathcal{R}$ is isomorphic to $K_p$.  
By Theorem~\ref{Gallai},  all edges in $\mathcal{R}$ are colored by at most two colors under $c$.  One can see that any monochromatic $H$ in $\mathcal{R}$ under $c$ will result in a monochromatic $H$ in $G$ under $c$. It is not   surprising   that  Gallai-Ramsey numbers $GR_k(H)$ are related to  the classical Ramsey numbers $R_2(H)$.  Recently,  Fox,  Grinshpun and  Pach posed the following  conjecture on $GR_k(H)$ when $H$ is a complete graph. 

\begin{conj}[\cite{FGP}]\label{Fox} For all integers $k\ge1$ and $t\ge3$,
\[
GR_k(K_t) = \begin{cases}
			(R_2(K_t)-1)^{k/2} + 1 & \text{if } k \text{ is even} \\
			(t-1)  (R_2(K_t)-1)^{(k-1)/2} + 1 & \text{if } k \text{ is odd.}
			\end{cases}
\]
\end{conj}

The first case of Conjecture~\ref{Fox} follows from a result of     Chung and Graham~\cite{chgr} in 1983.  The next open case when $t=4$ was recently settled in~\cite{K4}. 
In this paper, we study  Gallai-Ramsey numbers of  odd cycles. Using  the same construction given by   Erd\H{o}s, Faudree, Rousseau and Schelp  in 1976 (see Section 2  in \cite{EFRS}) for classical Ramsey numbers of odd cycles,   we see that $GR_k(C_{2n+1})\ge n \cdot 2^k + 1$  for all  $k \ge 1$ and $n \ge 2$.  
 General    upper bounds for   $GR_k (C_{2n+1})$ were first studied in \cite{c5c6} and later   improved  in  \cite{Hall}.   

\begin{thm}[\cite{Hall}]\label{general}
For all   $k \ge 1$ and $n \ge 2$, 
\[
 n \cdot 2^k + 1\le GR_k(C_{2n+1}) \le (2^{k+  3} -   3)n\,   \ln n.
\]
\end{thm}

Theorem~\ref{C3} and Theorem~\ref{C5}  below determine   the exact values of $GR_k(C_3)$ and $GR_k(C_5)$, respectively. A simpler proof of Theorem~\ref{C3} can be found in~\cite{exponential}.

\begin{thm}[\cite{chgr}]\label{C3} For all $k \ge 1$, 
	$GR_k(C_3) = \begin{cases}
			5^{k/2} + 1 & \text{if } k \text{ is even} \\
			2 \cdot  5^{(k-1)/2} + 1 & \text{if } k \text{ is odd.}
			\end{cases}$
\end{thm}

\begin{thm}[\cite{c5c6}]\label{C5}
  For all $k \ge 1$, $GR_k(C_{5}) = 2 \cdot 2^k + 1$.
\end{thm}

Recently, Bruce and Song~\cite{DylanSong} considered the next step and determined  the exact values of $GR_k(C_7)$ for all integers $k\ge1$. 
\begin{thm}[\cite{DylanSong}]\label{C7}
 For every  integer $k \ge 1$,   $GR_k(C_{7}) = 3 \cdot 2^k + 1$.
\end{thm}

We continue to study the Gallai-Ramsey numbers of odd cycles in this paper.  
We determine the exact values of Gallai-Ramsey numbers of $C_9$ and $C_{11}$ in this paper by showing that the lower bound in Theorem~\ref{general} is also the desired upper bound.  That is,  
we prove that  $GR_k(C_{2n+1}) \le  n \cdot 2^k + 1$ for all integers $n\in\{4,5\}$ and  $k \ge 1$.    Jointly with Bosse and Zhang~\cite{C13C15}, we are currently working on  the   Gallai-Ramsey numbers of $C_{13}$ and $C_{15}$, using   the key ideas developed in this paper. We believe the method we developed in this paper  and \cite{C13C15} will be helpful in determining the exact values of Gallai-Ramsey numbers of $C_{2n+1}$ for   all   $n\ge 8$.  Theorem~\ref{C9C11}  is our main result.   

\begin{thm}\label{C9C11}
For all integers $n\in\{4,5\}$ and  $k \ge 1$, $GR_k(C_{2n+1}) = n \cdot 2^k + 1$.
\end{thm}

It is worth mentioning that Theorem~\ref{C9C11}  also provides partial evidence for the first two open cases of the Triple Odd Cycle Conjecture due to Bondy and Erd\H{o}s~\cite{BE}, which states that $R_3(C_{2n+1})=8n+1$ for all integers  $n\ge2$. \L uczak~\cite{Luczak} showed that  $R_3(C_{2n+1}) =8n + o(n)$, as $n\rightarrow \infty$,  and Kohayakawa,
Simonovits and Skokan~\cite{TOCC} announced a proof in 2005 that the Triple Odd Cycle Conjecture holds when $n$ is sufficiently large.   

We shall make  use of the following   known results in the proof of Theorem~\ref{C9C11}.  

\begin{thm}[\cite{BE}]\label{2n+1} For all $n \ge 2$, 
 $R_2(C_{2n+1}) =4n+1$.
\end{thm}

\begin{prop}[\cite{CS}]\label{C4C6}
$R_2(C_{4}) =6$ and $R_2(C_{6}) =8$.
\end{prop}

Finally, we need to introduce more notation.   For positive integers $n, k$ and a complete graph $G$, let $c$ be any Gallai $k$-coloring of $G$ with color classes $E_1, \dots, E_k$. Then $c$ is  \dfn{bad}  if $G$ contains  no   monochromatic $C_{2n+1}$ under  $c$. 
For any  $W\subseteq V(G)$ and  any color $i\in [k]$,    $E:=E_i\cap E(G[W])$  is an \dfn{induced matching} in $G[W]$ if  $E$ is a matching in $G[W]$.
For  two disjoint sets $A, B\subseteq V(G)$,  $A$ is \dfn{mc-complete} to $B$ under the coloring $c$ if all the edges between $A$ and $B$  in $G$ are colored the same color under $c$;  and  we simply say $A$ is     \dfn{$j$-complete} to $B$    if all the edges between $A$ and $B$  in $G$ are colored by some color $j\in[k]$    under $c$;  and   $A$ is \dfn{blue-complete}     to $B$    if all the edges between $A$ and $B$  in $G$ are colored  blue  under $c$. 
For convenience, we use   $A \less B$ to denote  $A -B$; and $A \less b$ to denote  $A -\{b\}$ when $B=\{b\}$. We conclude this section with two  useful lemmas.

\begin{lem}\label{n,n+1 Lemma} For all integers $n\ge3$ and $k\ge1$, let $c$ be a   $k$-coloring of the edges of a complete graph $G$  on at least $2n+1$ vertices.   Let  $Y, Z\subseteq V(G)$ be   two disjoint sets with $|Y|  \ge   n$ and $ |Z| \ge n$. If   $Y$ is mc-complete, say blue-complete, to  $Z$  under the coloring $c$,  then no vertex  in $ V(G) \setminus (Y \cup Z)$ is blue-complete to $Y \cup Z$ in $G$.  Moreover, if $|Z| \ge n+1$, then $G[Z]$   has no blue edges. Similarly, if  $|Y| \ge n+1$, then   $G[Y]$ has no blue edges.
\end{lem}

\pf  Suppose  there exists a vertex $x\in  V(G) \setminus (Y \cup Z)$ such that $x$ is  blue-complete to $Y \cup Z$ in $G$.
Let $Y=\{y_1, \ldots, y_{_{|Y|}}\}$ and $Z=\{z_1, \ldots, z_{_{|Z|}}\}$. We may further assume that $z_1z_2$  is colored blue under $c$  if  $|Z| \ge n+1$ and $G[Z]$ has a blue edge.  We then obtain a blue $C_{2n+1}$ with  vertices $y_1, x, z_1, y_2,z_2, \ldots, y_n, z_n$ in order when $|Y|\ge n,|Z| \ge n$ or vertices  $y_1, z_1, z_2, y_2, z_3, \ldots, y_n, z_{n+1}$ in order when $|Z| \ge n+1$ and $G[Z]$ has a  blue edge $z_1z_2$, a contradiction. Thus  no vertex  in $ V(G) \setminus (Y \cup Z)$ is blue-complete to $Y \cup Z$ in $G$; and if  $|Z| \ge n+1$, then $G[Z]$   has no blue edges. Similarly, one can prove that if  $|Y| \ge n+1$, then   $G[Y]$ has no blue edges. \qed\medskip


\begin{lem}\label{Hamilton}
For all integers  $\ell\ge3$ and $n\ge1$,   let  $n_1, n_2, \ldots, n_\ell$ be positive integers such that $n_i\le n$ for all $i\in[\ell]$ and $n_1+n_2+\cdots+n_\ell\ge2n+1$. Then the complete multipartite graph $K_{n_1, n_2, \ldots, n_\ell}$ has a cycle of length $2n+1$. \end{lem} 
\begin{pf}
Let $G:=K_{n'_1, n'_2, \ldots, n'_\ell}$ be an induced   subgraph of $K_{n_1, n_2, \ldots, n_\ell}$ with $\ell'\ge3$, $n'_1+n'_2+\cdots+n'_\ell=2n+1$ and for all $i\in[\ell']$, $1\le n_i'\le n$. Then $\delta(G) \ge n+1 \ge |G|/2$.   By a well-known theorem of Dirac~\cite{Dirac}, $G$ has a Hamilton  cycle, and so $K_{n_1, n_2, \ldots, n_\ell}$ has a cycle of length $2n+1$.\qed\\
\end{pf}

\section{Proof of Theorem \ref{C9C11}}
Let $n\in\{4,5\}$. By  the construction given by   Erd\H{o}s, Faudree, Rousseau and Schelp  in 1976 (see Section 2  in \cite{EFRS}) for classical Ramsey numbers of odd cycles,   $GR_k(C_{2n+1}) \ge n \cdot 2^k + 1$ for    all $k \ge 1$. We next show that $GR_k(C_{2n+1}) \le n \cdot 2^k + 1$ for all $k \ge 1$.  This is trivially true for $k=1$. By Theorem \ref{2n+1} and the fact that $GR_2(C_{2n+1})=R_2(C_{2n+1})$, we may assume that $k\ge3$. Let $G:=K_{n\cdot 2^k + 1}$ and let $c$ be any Gallai $k$-coloring of $G$.  
We next show that  $G$ contains a  monochromatic copy of $C_{2n+1}$  under the coloring $c$. 
	
 Suppose that  $G$ does not  contain any  monochromatic $C_{2n+1}$ under $c$.  Then $c$ is  bad. Among all complete graphs on $n\cdot 2^k+1$ vertices with a bad Gallai $k$-coloring,  we choose $G$ with $k$ minimum.  We next  prove a series of  claims.	

\begin{claim}\label{l-vertex}
Let $W\subseteq V(G)$  and let $\ell\ge 3$ be an integer.  Let $x_1, \ldots, x_\ell \in V(G) \setminus W$ such that   $\{x_1, \ldots, x_\ell\}$ is mc-complete, say blue-complete, to  $W$  under $c$. Let    $q\in\{0, 1, \ldots, k-1\} $ be the number of colors, other than blue,  missing on $ G[W]$ under $c$.   
\begin{enumerate}[(i)]

\item If $\ell \ge n $, then $|W| \le  n \cdot 2^{k-1-q}$. 
 
\item If $\ell = n-1$, then $|W| \le  
n \cdot 2^{k-1-q} +2$.
\item If $\ell = n-2$, then  $n=5$ and $|W| \le 8 \cdot 2^{k-1-q} - 1$.
\end{enumerate}
\end{claim}

\pf    The statement in each of (i), (ii) and  (iii)  is trivially  true if  $|W| < \max\{2n+1 -\ell, n+1\}$. So we may assume that  $|W| \ge \max\{2n+1 -\ell, n+1\}$. We may further assume that  $G[W]$ contains at least one  blue edge, else,  by minimality of $k$, $|W|\le n \cdot 2^{k-1-q}$, giving the result.   
Note that $q\le k-1$. If $q=k-1$, then     all the edges of $G[W]$ are colored only blue. Since    $\{x_1,\ldots, x_\ell\}$ is blue-complete to $W$ and   $|W| \ge  \max\{2n+1 -\ell, n+1\}$, we see that   $G[W\cup \{x_1,  \ldots, x_\ell\}]$ contains a   blue $C_{2n+1}$, a contradiction.  Thus   $q \le k-2$.    
Since $|W|\ge n+1$ and $G[W]$ contains at least one  blue edge, by Lemma \ref{n,n+1 Lemma},  $\ell\le n-1$.  Let $W^*$ be a minimal set of vertices in $W$ such that      $G[W\less W^*]$ has no blue edges. By  minimality of $k$, $|W\less W^*|  \le n \cdot 2^{k-1-q}$. 

 We now consider the case when $\ell = n-1$. Then $|W|\ge 2n+1-\ell=n+2$.  
If $G[W]$   contains three blue edges, say  $u_1v_1, u_2v_2, u_3v_3$,  such that $u_1, u_2, u_3, v_1, v_2, v_3$ are all distinct, then   we obtain a blue $C_{2n+1}$ with vertices $x_1, u_1, v_1, x_2$, $u_2, v_2, x_3, u_3, v_3$ in order  (when $n=4$) and vertices $x_1, u_1, v_1, x_2, u_2, v_2, x_3, u_3, v_3, x_4, u $ in order (when $n=5$,  where $u\in W\less \{u_1, u_2, u_3, v_1, v_2, v_3\}$), a contradiction. Thus   $|W^*|\le 2$, and so   $|W|\le n \cdot 2^{k-1-q} + 2$.  

\begin{figure}[h]
\centering
\includegraphics[scale=1.2]{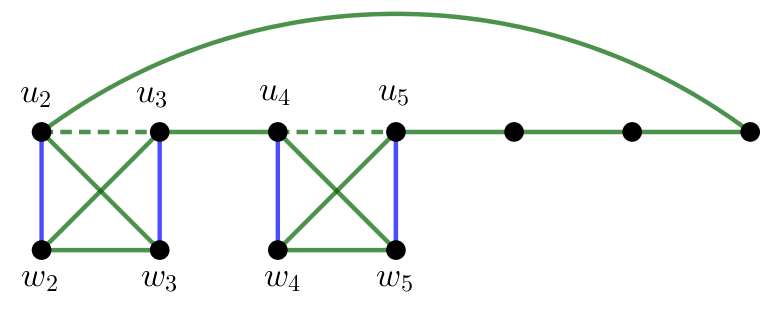}
\caption{An example of a  green $C_{11}$ arising  from the  green $C_7$.}
\label{C7toC11}
\end{figure}

It remains to  consider the case when $3\le \ell\le n-2$. Then   $n=5$ and $\ell=n-2=3$. Note that $|W|\ge 2n+1-\ell\ge8$. 
Let $P$ be a longest blue path in $G[W]$  with vertices $v_1, \dots, v_{_{|P|}}$ in order.   Since  $\{x_1,x_2, x_{3}\}$ is blue-complete to $W$, we see that  $|P|\le 5$, else we obtain a blue $C_{11}$ with vertices  $x_1, v_1, \ldots, v_{6}, x_2, u_1, x_{3}, u_{2} $ in order, where $u_1,   u_{2}\in W\less \{v_1, \ldots, v_{6}\}$, a contradiction. 
Assume first that  $|W^*|\le 4$. Then,
\[
|W| = |W \less W^*| + |W^*| \le n \cdot 2^{k-1-q} + 4 < 8 \cdot 2^{k-1-q} - 1,
\]
because $q \le k-2$ and $k\ge3$.
So we may assume that $|W^*|\ge5$. By the choice of $W^*$,  we see that  $|P|\in\{2,3\}$, else we obtain a blue $C_{11}$. Furthermore, if $|P|=3$, then $G[W\less V(P)]$ has no blue path on three vertices. Thus all the blue edges in $G[W\less V(P)]$ induce a blue matching. Let $m:=|W^*\less V(P)|$ and 
 let $u_2w_2, \ldots, u_{m+1}w_{m+1}$   be    all the   blue edges  in $G[W\less V(P)]$, where $u_2, \ldots, u_{m+1}, w_2, \ldots, w_{m+1}$ are all distinct. By the choice of $W^*$, we may assume that $u_2, \ldots, u_{m+1}\in W^*$. Let $u_1=v_1$ and $w_1=v_2$, and  $A:=W\less (V(P)\cup \{u_2, \ldots, u_{m+1}, w_2, \ldots, w_{m+1}\})$.  Let $B:=\{u_1, u_2, \ldots, u_{m+1}\}$ when $| A|\le1$ and let $B:=\{u_1, u_2, \ldots, u_{m+1}\}\cup\{a_1, a_2\}$ when $|A|\ge2$ and $a_1, a_2\in A$ with $a_1\ne a_2$.  We claim that   $|B|\le 3 \cdot 2^{k-1-q}$. Suppose $|B|\ge 3 \cdot 2^{k-1-q}+1$.  By   Theorem~\ref{C7}, $G[B]$ has a monochromatic, say green,  $C_7$.  Then $|V(C_7)\cap \{u_1, u_2, \ldots, u_{m+1}\}|\ge5$ and so $C_7\less \{a_1, a_2\}$ has a matching of size two. We may assume that $u_2u_3, u_4u_5\in E(C_7)$. Since $G$ has no rainbow triangles under the coloring $c$, we see that for any $i\in\{2,4\}$, $\{u_i, w_i\}$ is green-complete to $\{u_{i+1}, w_{i+1}\}$. Thus we obtain a green $C_{11}$ from the $C_7$ by replacing the edge $u_2u_3$ with the path $u_2w_3w_2 u_3$ and edge $u_4u_5$ with the path $u_4w_5 w_4 u_5$, a contradiction (see Figure~\ref{C7toC11}).  Thus $|B|\le 3 \cdot 2^{k-1-q}$, as claimed. \\

 When   $|A|\le1$, we have  $|W|= |A|+2|B|+|V(P)\less\{ v_1,v_2 \}|\le 1+6 \cdot 2^{k-1-q}+1<8 \cdot 2^{k-1-q} - 1$ because $q\le k-2$ and $k\ge3$. When $|A|\ge2$,  
since $G[A\cup\{w_1,  w_2, \ldots, w_{m+1}\}]$ has no blue edges, by minimality of $k$, 
$|A\cup \{w_1,  w_2, \ldots, w_{m+1}\}|\le 5\cdot 2^{k-1-q}$. Hence, 
\begin{align*}
|W|&= |A\cup\{w_1,  w_2, \ldots, w_{m+1}\}|+|B\less \{a_1, a_2\}|+|V(P)\less\{ v_1,v_2 \}|\\
&\le 5 \cdot 2^{k-1-q}+(3 \cdot 2^{k-1-q}-2)+1\\
&= 8 \cdot 2^{k-1-q} - 1. 
\end{align*}
This completes the proof of Claim \ref{l-vertex}.  \qed\\

Let $X_1,   \ldots, X_m$  be a maximum sequence of disjoint subsets of $V(G)$ such that,  for all $j \in [m]$, one of the following holds. 
\begin{enumerate}[(a)]
\item    $1 \le |X_j| \le 2$, and   $X_j$ is mc-complete to  $V(G) \setminus \bigcup_{i \in [j]} X_i$ under $c$, or

\item     $3 \le |X_j| \le 4$,  and $X_j$ can be partitioned into two non-empty sets $X_{j_1}$ and $X_{j_2}$, where $j_1, j_2\in [k]$ are two distinct colors, such that for each  $t\in\{1,2\}$,   $1\le |X_{j_t}|  \le 2$,    $X_{j_t}$ is $j_t$-complete to $V(G) \setminus \bigcup_{i \in [j]} X_i$ but  not $j_t$-complete to $X_{j_{3-t}}$, and all the edges between  $X_{j_1}$ and  $X_{j_2}$ in $G$ are colored using only the colors $j_1$ and $j_2$.  
\end{enumerate}

  \noindent Note that such a sequence $X_1,   \ldots, X_m$  may not exist. 
  Let $X:= \bigcup_{j \in [m]} X_j$.   For each $x\in X$,  let $c(x)$ be the unique color on the edges between $x$ and $V(G) \less X$ under $c$.  For all $i \in [k]$, let $X_i^*:= \{x \in X : c(x) = \text{color } i\}$. Then  $X = \bigcup_{i \in [k]} X_i^*$. It is worth noting that for all $i \in [k]$, $X_i^*$ is possibly empty.  By abusing the notation, we use $X_b^*$ to denote $X_i^*$  when the color $i$ is blue. Similarly, we use $X_r^*$ to denote  $X_i^*$  when the color $i$ is red.

\begin{claim}\label{X*}
For all $i \in [k]$, $|X_i^*| \le 2$.
\end{claim}

\begin{pf}
Suppose the statement is false.  Then $m \ge 2$.  When choosing $X_1, X_2, \ldots, X_m$, let $j\in[m-1]$ be the largest index such that $|X_p^*\cap (X_1\cup X_2\cup\dots \cup X_{j})|\le2$ for all $p\in[k]$. Then     $3\le |X_i^*\cap (X_1\cup X_2\cup\dots \cup X_j\cup X_{j+1})|\le4$ for some color $i\in[k]$ by the choice of $j$.  Such a color $i$ and an index $j$ exist due to the assumption that the statement of Claim \ref{X*} is false.  Let $A :=X_1\cup X_2\cup  \cdots\cup X_j\cup X_{j+1}$. 
By the choice of   $X_1, X_2, \ldots, X_m$,  there are   at most two colors $i\in[k]$ such that $3\le |X_i^*\cap A|\le4$. We may assume that such a color $i$ is either blue or red. 
Let  $A_b:=\{x\in A: \, c(x) \text{ is color blue} \}$ and $A_r:=\{x\in A: \, c(x) \text{ is color red} \}$.  It suffices to consider the worst case when    $3\le |A_b|  \le4$ and   $3\le   |A_r| \le4$. Then    for any color $p\in[k]$ other than red and blue,  $|X_p^*\cap A|\le 2$.  Thus by the choice of $j$,  $|A\less (A_b \cup A_r)|\le 2(k-2)$. We may assume that $|A_b|\ge|A_r|$. Note that $|A_b|\le n$.  If $|A_b|\ge n-1$, then by Claim \ref{l-vertex}(ii) applied to any   $n-1$ vertices in $A_b$ and  $V(G)\less A$,  we see that $|V(G) \setminus A|\le n \cdot 2^{k-1} +2 $. Thus, 
\[
|G| = |A\less (A_b \cup A_r)| + |A_b| + |A_r|+ |V(G) \setminus A| \le 2(k-2) + n+n+ (n \cdot 2^{k-1} +2) < n \cdot 2^k + 1
\]
for all $k \ge 3$ and  $n\in\{4,5\}$, a contradiction. Thus   $ 3\le |A_b| \le n-2$.  Then $|A_b|=3$ and $n=5$.  By Claim \ref{l-vertex}(iii) applied to $A_b$ and  $V(G)\less A$, we see that $|V(G) \setminus A|\le 8 \cdot 2^{k-1} -1 $. Thus, 
\[
|G| = |A\less (A_b \cup A_r)| + |A_b| + |A_r|+ |V(G) \setminus A| \le 2(k-2) + 3+3+ (8 \cdot 2^{k-1} -1) < 5 \cdot 2^k + 1
\]
for all $k \ge 3$, a contradiction.  \qed\\
\end{pf}

By Claim \ref{X*}, $|X| \le 2k$.  Let $X'\subseteq X$  be such that for  all $i \in [k]$,  $|X'\cap X^*_i|=1$ when $X^*_i \ne \es$.  Let $X'':= X\less X'$. 
Now consider a  Gallai partition $A_1, \ldots, A_p$ of $G \setminus X$ with $p\ge2$. 
We may  assume that  $1\le |A_1| \le \cdots \le  |A_s|<3 \le  |A_{s+1}| \leq \cdots \leq |A_{p}|$, where $0\le  s\le p$. Let  $\mathcal{R}$ be  the reduced graph of $G\less X$ with vertices $a_1, a_2, \dots, a_{p}$, where $a_i\in A_i$ for all $i\in[p]$. By Theorem~\ref{Gallai}, we may assume that  the edges of $\mathcal{R}$ are colored red and blue. Note that any monochromatic $C_{2n+1}$ in $\mathcal{R}$ would yield a monochromatic $C_{2n+1}$ in $G$. Thus $\mathcal{R}$ has neither a red nor a blue $C_{2n+1}$.  By Theorem~\ref{2n+1}, $p \le 4n$. Then $|A_p|\ge 2$ because $|G\less X| \ge n\cdot 2^k+1 -2k\ge 8n-5$. If $|A_p|= 2$, then $k=3$.  Thus $|A_{p-4n+8}|=2$, else $ |G|\le 2(4n-8)+(p-(4n-8))+|X|\le 8n-2<n\cdot 2^3+1$, a contradiction. Since $R_2(C_{2n-3})=4n-7$ by Theorem~\ref{2n+1}, we see that $\mc{R}[\{a_{p-4n+8}, a_{p-4n+9}, \ldots, a_p\}]$ has a monochromatic, say blue,  $C_{2n-3}$, and so $G[A_{p-4n+8}\cup A_{p-4n+9}\cup \cdots\cup A_p]$ has a blue $C_{2n+1}$, a contradiction.  Thus $|A_p|\ge3$ and so   $p-s\ge1$. Let 
\[
\begin{split}
B &:=  \{a_i \in \{a_1, \ldots, a_{p-1}\} \mid a_ia_1 \text{ is colored blue in } \mc{R} \}\\
R &:= \{a_j \in \{a_1, \ldots, a_{p-1}\} \mid a_ja_1 \text{ is colored red in } \mc{R} \}
\end{split}
\]
Then $|B|+|R|=p-1$. 
Let $B_G:= \bigcup_{a_i \in B} A_i$ and  $R_G:=\bigcup_{a_j\in R} A_j$. 

\begin{claim}\label{B}
If $|A_p|\ge n$ and $|B|\ge 3$ (resp. $|R|\ge 3$), then $|B_G|\le 2n$ (resp. $|R_G|\le 2n$). 
\end{claim}

\begin{pf}  
Suppose $|A_p|\ge n$ and $|B|\ge 3$ but $|B_G|\ge 2n+1$. By Lemma \ref{n,n+1 Lemma},  $G[B_G]$ has no blue edges and no vertex  in $ X$ is blue-complete to $V(G)\less X$.  Thus all the edges of $\mathcal{R}[B]$ are colored red in $\mathcal {R}$. Let $m:=|B| $   and let $B:=\{a_{i_1}, a_{i_2}, \ldots, a_{i_m}\}$ with $|A_{i_1}|\ge |A_{i_2}|\ge\cdots\ge |A_{i_m}|$. Then $G[B_G]- \bigcup_{j=1}^{m}E(G[A_{i_j}])$ is a complete multipartite graph with at least three parts.  If $|A_{i_1}|\le n$, then by Lemma~\ref{Hamilton} applied to $G[B_G]- \bigcup_{j=1}^{m}E(G[A_{i_j}])$, $G[B_G]$ has a red $C_{2n+1}$, a contradiction. 
 Thus  $|A_{i_1}|\ge n+1$.     Let $Q_b:=\{v\in R_G: v \text{ is blue-complete to } A_{i_1}\}$, and $Q_r:=\{v\in R_G: v \text{ is red-complete to } A_{i_1}\}$. Then $Q_b\cup Q_r=R_G$.  
 Let $Q: =(B_G\less A_{i_1})\cup Q_r \cup X^*_r$. Then $A_{i_1}$ is red-complete to $Q$ and $G[Q]$ must contain red edges,  because   $|B|\ge 3$ and all the edges of $\mathcal{R}[B]$ are colored red. By Lemma \ref{n,n+1 Lemma} applied to $A_{i_1}$ and $Q$, $|Q|\le n$. 
 Note that  $|A_p\cup Q_b|\ge|A_p|\ge|A_{i_1}|\ge n+1$ and $A_p\cup Q_b$ is blue-complete to $A_{i_1}$. By Lemma \ref{n,n+1 Lemma} applied to $A_{i_1}$ and $A_p\cup Q_b$,    $G[A_p\cup Q_b]$ has   no  blue edges. Since no vertex  in $ X$ is blue-complete to $V(G)\less X$, we see that  neither $G[A_p\cup Q_b\cup (X''\less X^*_r)]$ nor $G[B_G\cup X']$ (and thus $G[A_{i_1}\cup (X'\less X^*_r)]$)  has  blue edges. By minimality of $k$,   $|A_p\cup Q_b\cup (X''\less X^*_r)|\le n\cdot 2^{k-1}$.
 Suppose first that  $Q_r \cup X^*_r=\emptyset$. 
 Then $Q_b=R_G$, so that
\[ |G| =|B_G\cup X'|+|A_p\cup Q_b\cup X''|  \le  n\cdot 2^{k-1}+n\cdot 2^{k-1} <n\cdot 2^{k}+1,\]
 a contradiction. Thus $Q_r \cup X^*_r\ne\emptyset$. Since $|B|\ge3$, we see that $|B_G\less A_{i_1}|\ge2$. Thus  $n\ge |Q|\ge3$.  Since $G[A_{i_1}\cup (X'\less X^*_r)]$ has no blue edges, 
 by Claim \ref{l-vertex}  applied to $Q$ and  $ A_{i_1}$ we see that 
  \[
   |A_{i_1}\cup (X'\less X^*_r)|\le\begin{cases}  
 n \cdot 2^{k-2} +2, &  \text{if } |Q| \in\{ n-1, n\}\\
 8 \cdot 2^{k-2}-1, &  \text{if } |Q| = n-2 \text{ and } n=5.
 \end{cases}
 \] 
But then 
\begin{align*}
|G|&=|Q|+|A_{i_1}\cup (X'\less X^*_r)|+|A_p\cup Q_b\cup (X''\less X^*_r)|\\
  &\le\begin{cases}  
  n+(n \cdot 2^{k-2} +2) +n\cdot 2^{k-1}, &  \text{if } |Q| \in\{ n-1, n\}\\
3+ (8 \cdot 2^{k-2}-1)+n\cdot 2^{k-1}, &  \text{if } |Q| = n-2 \text{ and } n=5.
 \end{cases}\\
 &<n\cdot 2^{k}+1
\end{align*}
for all $k\ge3$ and $n\in\{4,5\}$,  a contradiction. Hence, $|B_G|\le 2n$. Similarly, one can prove that if $|A_p|\ge n$ and $|R|\ge 3$, then $|R_G|\le 2n$. \qed
\end{pf}

\begin{claim}\label{n=5, p}
$p\le 2n-1$.
\end{claim}
\begin{pf}
Suppose $p \ge 2n$.  Then $|B| + |R| = p-1\ge 2n-1$. We claim that  $|A_p| \le n-1$. Suppose $|A_p| \ge n$.  We may assume that $|B| \ge |R|$.  Then    $|B_G| \ge |B| \ge n>3$.  
By Claim \ref{B}, $|B_G| \le 2n$.  If $|R_G|\ge n+1$, then by Lemma \ref{n,n+1 Lemma} to $A_p$ and $R_G$, $G[R_G]$ has no red edges, and no vertex in $X$ is red-complete to $V(G) \less X$. Then $|X''|\le k-1$ and $G[R_G \cup X']$ has no red edges. By minimality of $k$,  $|R_G \cup X'| \le n \cdot 2^{k-1}$.  Then  
\[
|A_p|=|G|-|B_G|-|R_G \cup X'|-|X''|\ge n\cdot 2^k+1-2n-n \cdot 2^{k-1}-(k-1)\ge 2n-1,\]
  for all $k\ge3$. By Lemma~\ref{n,n+1 Lemma} applied to $A_p$ and $B_G$,   $G[A_p]$ has no blue edges and no vertex in $ X$ is blue-complete to $V(G)\less X$.   Thus   $G[A_p\cup X'']$ has neither red nor blue edges, and so      $|A_p\cup X''|\le n\cdot 2^{k-2}$ by the choice of $k$.   But then  
  \[
  |B_G|= |G|-|R_G \cup X'|-|A_p\cup X''|\ge n\cdot 2^k+1-n \cdot 2^{k-1}-n\cdot 2^{k-2} \ge 2n + 1,
  \]
   contrary to Claim \ref{B}. This proves that 
$|R_G|\le n$. Then 
\[
|A_p\cup X'|=|G|-|B_G|-|R_G|-|X''|\ge (n\cdot 2^k+1)-2n -n-k>n\cdot 2^{k-1}+1.
\]
By minimality of $k$,  $G[A_p\cup X']$ must have  blue edges. Since $|A_p|\ge n$ and $|B_G|\ge n$, by Lemma~\ref{n,n+1 Lemma} applied to $A_p$ and $B_G$, $|A_p|=n$ and no vertex in $X$ is blue-complete to $V(G)\less X$. Thus   $|X|\le 2(k-1)$. But then 
\[
|G|=|B_G|+|R_G|+|A_p|+|X|\le 2n+n+n+2(k-1) < n\cdot 2^k+1,\]
for all $k\ge3$, a contradiction.  
  This proves that 
 $|A_p| \le n-1$, as claimed.  
 
  Since $ |A_p| \ge3$, we have  $3 \le |A_p| \le n-1$.  Then   $k = 3$ because $n\in\{4,5\}$ and $|G|=n\cdot 2^k+1$.  It follows that   $|G| = 8n+1$ and $|X| \le 6$. Therefore, $|B_G| + |R_G| = |G| - |A_p| - |X| \ge (8n+1) - (n-1)-6 =7n-4$.  We may thus assume that $|B_G| > 2n+3$. We next prove that $|A_p|\le n-2$. Suppose $|A_p| = n-1$.  If $G[B_G]$   contains three blue edges $u_1v_1, u_2v_2, u_3v_3$ such that $u_1, u_2, u_3, v_1, v_2, v_3$ are all distinct, then   we obtain a blue $C_{2n+1}$ with vertices in  $A_p\cup   \{u_1, u_2, u_3, u_4, v_1, v_2, v_3\}$, where $u_4\in B_G\less   \{u_1, u_2, u_3,  v_1, v_2, v_3\}$, a contradiction. Thus there exists $B^*\subseteq B_G$ such that $|B^*|\le2$ and $G[B_G\less B^*]$ has no blue edges.
  Then  $|B_G \less B^*| > 2n+1$, and so  $|B\less B^*|\ge3$ because $|A_i|\le n-1$ for all $i\in[p]$. By the choice of $B^*$,   all the edges in $\mc{R}[B\less B^*]$ are colored red.   But then by Lemma \ref{Hamilton}, $G[B_G \less B^*]$ has  a red $C_{2n+1}$, a contradiction.
This proves that  $3\le |A_p| \le n-2$.  Then $|A_p|=3$, $ n=5$,     $|G|=41$, and $p\le 20$. 
 If    $|A_{p-7}|=3$ or  $|A_{p-12}|\ge2$, then  $\mc{R}[\{a_{p-8}, a_{p-7}, \ldots,    a_{p}\}]$ has a monochromatic   $C_5$,  or  $\mc{R}[\{ a_{p-12},a_{p-11}, \ldots,   a_{p}\}]$ has a monochromatic  $C_{7}$ because $R_2(C_5)=9$ and $R_2(C_7)=13$.  In either case, we see that  $G$ has  a monochromatic $C_{11}$, a contradiction.  Thus  $|A_{p-7}|\le 2$ and $|A_{p-12}|\le1$.  
Then  $|A_{p-7}| = 2$, else $|G| \le 7\cdot 3 + 13 \cdot 1 + 6 < 41$, a contradiction. 
Since $R_2(C_6)=8$,  we see that  $\mc{R}[\{a_{p-7}, a_{p-6}, \ldots, a_{p}\}]$  has a monochromatic, say blue, $C_{6}$, and so $G\less X$ has a blue $C_{10}$. Thus no vertex in $X$ is blue-complete to $G\less X$ and so $|X|\le 2(k-1)=4$. 
Furthermore, if $|A_{p-8}|=2$, then 
$|A_{p-4}| = 2$, else $\mc{R}[\{a_{p-8}, a_{p-7}, \ldots,    a_{p}\}]$ has a monochromatic   $C_5$, and so $G$ has  a monochromatic $C_{11}$, a contradiction.  But then 
 $|G| \le 4\cdot 3 + 8\cdot 2+ (p-12) \cdot 1 + |X| \le 40< 41$ when $|A_{p-8}|=2$; and  $|G| \le 7\cdot 3 +   2+ (p-8) \cdot 1 + |X| \le 39< 41$ when $|A_{p-8}|\le1$. In both cases, we obtain  a contradiction. 
  \qed 
\end{pf}

\begin{claim}\label{Ap}
$|A_p|\ge n+1$.
\end{claim}

\pf  Suppose $|A_p|\le n$.  By Claim \ref{n=5, p}, $p \le 2n-1$.  We may assume that $a_pa_{p-1}$ is colored blue in $\mc{R}$. Then $|A_p\cup A_{p-1}\cup X|\le 2n+2(k-1)$, else we obtain a blue $C_{2n+1}$. 
If  $|A_{p-4}|\ge n-1$, then     $\mc{R}[\{a_{p-4}, a_{p-3}, \ldots,  a_{p}\}]$ has a monochromatic $C_3$ or $C_5$, and so   $G$ contains   a monochromatic  $C_{2n+1}$, a contradiction. Thus $|A_{p-4}|\le n-2$.   But then 
\[
 |G|\le (2n+2(k-1))+2n+(p-4)(n-2) \le 4n+(2n-5)(n-2) +2k-2<n\cdot 2^k+1.
\]
 for all $n\in\{4,5\}$ and $k\ge3$, a contradiction. 
  \qed \\

For the remainder of the proof, let  $B_G^*:=B_G\cup X_b^*$   and   $R_G^*:=R_G\cup X_r^*$.

\begin{claim}\label{one}
$2\le p-s\le 3n-7$.
\end{claim}

\pf  Suppose   $p-s\ge3n-6$. Then  $\mc{R}[\{a_{p-3n+7}, a_{p-3n+8}, \ldots,  a_{p}\}]$ has a  monochromatic    $C_{2n-5}$ because $R_2(C_{2n-5})=3n-6$ when $n\in\{4,5\}$. But then $G$   would  contain a monochromatic  $C_{2n+1}$. \medskip

Next suppose $p-s\le1$. Then $p-s=1$ because $p-s\ge1$.  Thus $|A_i|\le2$ for all $i\in[p-1]$ by the choice of $p$ and $s$.  By Claim \ref{n=5, p}, $p \le 2n-1$.    Then $|B_G \cup R_G|\le 2(p-1)$ and so $|B^*_G \cup R^*_G| \le 2(p-1)+2+2=2(p+1)\le 4n$. We may assume that $|B^*_G| \ge |R^*_G|$.  If $|R^*_G|\ge n$, then $|B^*_G|\ge n$. By Claim \ref{Ap} and Lemma \ref{n,n+1 Lemma}, $G[A_p]$ has neither  blue nor red edges. By minimality of $k$, $|A_p  | \le n \cdot 2^{k-2}$. But then 
\[
|G|=|B^*_G \cup R^*_G|+|A_p |+ |X \less (B^*_G \cup R_G^*)|\le 4n + n \cdot 2^{k-2} +2(k-2)< n\cdot 2^k+1
\]
for all $k\ge3$, a contradiction.  Thus $|R^*_G|\le n-1$. We claim that   $|B^*_G|\le 2n+2$. This is trivially true if $|B|\le n$.   If $|B|\ge n+1$, then  $|B_G|\le 2n$  by Claim~\ref{B}.  Thus $|B^*_G|\le 2n+2$, as claimed.   If $|B_G^*| \ge n-1$, then applying Claim \ref{l-vertex}(i,ii) to   $B_G^*$ and $A_p$ implies that 
 \[
  |B_G^*| + |A_p|\le\begin{cases}  
 (n-1)+(n \cdot 2^{k-1} +2), &  \text{if } |B_G^*| = n-1\\
  (2n+2)+n \cdot 2^{k-1}, &  \text{if } |B_G^*| \ge n.
 \end{cases}
 \] 
In either case, $ |B_G^*| + |A_p|\le2n+n \cdot 2^{k-1}+2$. But then 
\[
|G| = |R^*_G| + |B_G^*| + |A_p| + |X \less (B^*_G \cup R_G^*)| \le  (n-1) + (2n+n \cdot 2^{k-1} +2) +  2(k-2) <n \cdot 2^k + 1,
\]
for all $k\ge3$  and $n\in\{4,5\}$, a contradiction. 
Thus  $n-2\ge |B_G^*|\ge|R_G^*| $.  If   $|B_G^*| = 3$,  then $n=5$.  By Claim \ref{l-vertex}(iii) applied to  $B_G^*$ and $A_p$,  $|A_p| \le 8 \cdot 2^{k-1} -1$.  But then,
\[
|G| = |B^*_G| + |R_G^*| + |A_p| + |X\less (B^*_G \cup R_G^*)| \le 3 + 3 + (8 \cdot 2^{k-1} -1) +  2(k-2) < 5 \cdot 2^k + 1
\]
for all $k \ge 3$, a contradiction.  Thus $2\ge |B_G^*|\ge|R_G^*| $. Since $p\ge2$, we see that $B\ne\emptyset$ or $R\ne\emptyset$. 
Then by maximality of $m$ (see condition (a) when choosing $X_1, X_2, \ldots, X_m$), $B^*\ne \emptyset$, $R^*\ne \emptyset$, and $B^*_G$ is neither blue- nor red-complete to $R^*_G$ in $G$. But then, by maximality of $m$ again (see condition (b) when choosing $X_1, X_2, \ldots, X_m$), $B^*_G =\emptyset$ and   $R^*_G=\emptyset$,  contrary to   $p\ge2$.\qed 

\begin{claim}\label{A_{p-2}}
$|A_{p-2}|\le n-1$.
\end{claim}

\begin{pf}  Suppose $|A_{p-2}|\ge n$. Then $n \le |A_{p-2}| \le |A_{p-1}| \le |A_p|$ and so  $\mc{R}[\{a_{p-2}, a_{p-1},a_p\}]$ is not a monochromatic triangle in $\mc{R}$ (else we obtain a monochromatic $C_{2n+1}$).   Let $B_1$, $B_2$, $B_3$ be a permutation of $A_{p-2}$, $A_{p-1}$, $A_p$ such that $B_2$ is, say blue-complete,  to $B_1 \cup B_3$ in $G$. Then $B_1$ must be  red-complete to $B_3$ in $G$. We may assume that  $|B_1|\ge|B_3|$.  
 By Lemma~\ref{n,n+1 Lemma},    no vertex in $X$ is blue- or red-complete to $V(G)\less X$.     Let $A:=V(G)\less (B_1\cup B_2\cup B_3\cup X)$. Then by Lemma~\ref{n,n+1 Lemma},  no vertex in $A$ is red-complete to $B_1\cup B_3$ in $G$,  and no vertex in $A$ is blue-complete to $B_1\cup B_2$ or $B_2\cup B_3$ in $G$. This implies that 
$A$ must be  red-complete to $B_2$ in $G$. We next show that  $G[A]$ has no blue edges. 
Suppose that $G[A]$ has a blue edge, say, $uv$.  Let 
\[
\begin{split}
B_1^*&:=\{b\in A \mid  b \text{ is blue-complete to   }  B_1 \text{ only in } G \}\\
B_2^*&:=\{b\in A\mid  b \text{ is blue-complete to both   }  B_1 \text{ and }  B_3 \text{   in } G\}\\ 
B_3^*&:=\{b\in A\mid  b \text{ is blue-complete to   }  B_3 \text{ only in } G \}.\\ 
\end{split}
\]
Then $A=B_1^*\cup B_2^*\cup B_3^* $.  Note that $B_1^*,  B_2^*, B_3^* $ are  pairwise disjoint and possibly empty.   
 Let $b_1, \ldots, b_{n-1}\in B_1$, $ b_{n}, \ldots, b_{2n-2}\in B_2$, and $b_{2n-1}\in B_3$.  If $uv$ is an edge in $G[B_1^*\cup B_2^*]$, then we obtain a blue $C_{2n+1}$ with vertices $b_1, u, v, b_2, b_n, b_{2n-1}, b_{n+1}, b_3, b_{n+2}, \ldots, b_{n-1}, b_{2n-2}$
   in order, a contradiction. 
Similarly, $uv$ is not an edge in $G[B_2^*\cup B_3^*]$.  Thus $uv$ must be  an edge in $G[B_1^*\cup B_3^*]$ with one end in $B_1^*$ and the other in $B_3^*$.  We may assume that $u\in B_1^*$ and $v\in B_3^*$. Then we obtain a blue $C_{2n+1}$ with vertices $b_1, u, v, b_{2n-1}, b_n, b_2, b_{n+1},\ldots,  b_{n-1}, b_{2n-2}$ in order, a contradiction. This proves that $G[A]$ has no blue edges. By minimality of $k$, $|A| \le n\cdot 2^{k-1}$. 

We next show that $|B_2\cup A\cup X'|\le n \cdot 2^{k-1}$. Suppose $|B_2\cup A\cup X'|\ge n \cdot 2^{k-1}+1$. Then by minimality of $k$, $G[B_2\cup A\cup X']$ must contain blue edges. Since $G[A]$ has no blue edges, $A$ is red-complete to $B_2$,  and  no vertex in $X$ is blue-complete to $V(G)\less X$, we see that  $G[B_2]$ must contain blue edges. By Lemma~\ref{n,n+1 Lemma}, $|B_2|=n$. Then $B_2\ne A_p$. We may assume that $B_1=A_p$. By Lemma~\ref{n,n+1 Lemma}, $G[B_1]$ has neither blue nor red edges and so $G[B_1\cup X']$ has neither blue nor red edges. By minimality of $k$, $|B_1\cup X'|\le n\cdot 2^{k-2}$ and so $|B_3\cup X''|\le |B_1\cup X'|\le n\cdot 2^{k-2}$. Note that   $A=\emptyset$, else, let $v\in A$. Then $G[B_2\cup\{v\}]$ has blue edges and $B_2\cup\{v\}$ is blue-complete to either $B_1$ or $B_3$, contrary to Lemma~\ref{n,n+1 Lemma}.  But then 
\[
 |G|=|B_1\cup   X'|+|B_2|+|B_3\cup X''| \le n \cdot 2^{k-2}+n +n\cdot 2^{k-2} <n\cdot 2^{k}+1,
 \] 
 for all $k\ge3$, a contradiction.  
This proves that  $|B_2\cup A\cup X'|\le n\cdot 2^{k-1}$. 

 Since $|B_1|\ge |B_3|$ and  $|B_1|+|B_3|=|G|-|B_2\cup A\cup X'|-|X''|\ge n\cdot 2^{k-1}+1-(k-2)\ge 2n+1$, we see that $|B_1|\ge n+1$. Note that  $|B_2| \ge n$ and  $|B_3| \ge n$. By Lemma~\ref{n,n+1 Lemma},    $G[B_1]$  has neither red nor blue edges. Since each  vertex in $X$ is neither red- nor blue-complete to $B_1$,  $G[B_1\cup X'']$ has neither red nor blue edges. By minimality of $k$, $|B_1\cup X''|\le n\cdot 2^{k-2}$ and so $|B_3|\le |B_1|\le n \cdot 2^{k-2}$. But then  
 \[
 |G|=|B_2\cup A\cup X'|+|B_1\cup X''|+|B_3|\le n \cdot 2^{k-1}+n\cdot 2^{k-2}+n\cdot 2^{k-2}=n\cdot 2^{k},
 \] 
   a contradiction.  \qed \\
   \end{pf}

  By Claim~\ref{one},   $2 \le p-s \le 3n-7$ and so $|A_{p-1}|\ge3$.  We may now assume that $a_pa_{p-1}$ is colored blue in $\mc{R}$. Then   $a_{p-1}\in B$ and so $A_{p-1}\subseteq B_G$. Thus $|B_G|\ge|A_{p-1}|\ge3$.  

\begin{claim}\label{R*_G}
$|R^*_G|\le 2n$.
\end{claim}

\begin{pf} 
Suppose   $|R^*_G| \ge 2n+1$.  By Claim~\ref{Ap}, $|A_p| \ge n+1$. By Lemma~\ref{n,n+1 Lemma},  $G[R^*_G]$ has no red edges.  Thus $|R^*_G|=|R_G|$ and so no vertex in $X$ is red-complete to $V(G)\less X$. In particular, all the edges in $\mc{R}[R]$ are colored blue. By  Claim~\ref{B}, $|R| \le 2$.  By Claim~\ref{A_{p-2}},  $|A_{p-2}| \le n-1$.    Since $A_{p-1}\cap R_G=\emptyset$ and $|R_G| \ge 2n+1$, we see that $|R| \ge 3$, a contradiction.
\qed
\end{pf}

\begin{claim}\label{A_{p-1}}
$|A_{p-1}|\le n$.
\end{claim}
\begin{pf} 
Suppose   $|A_{p-1}| \ge n+1$.  
Then $|B_G| \ge |A_{p-1}|\ge n+1$.  By Lemma~\ref{n,n+1 Lemma},  neither $G[A_p]$ nor $G[B_G]$ has blue edges, and no vertex in $X$ is blue-complete to $V(G)\less X$.  Thus $|X|\le 2(k-1)$.  By the choice of $k$, $|B_G\cup X''|\le n \cdot 2^{k-1}$ and $|A_p\cup X'|\le n \cdot 2^{k-1}$. We claim that $G[R_G]$ has blue edges. Suppose  $G[R_G]$ has no blue edges. Then $G[A_p\cup R_G\cup X']$ has no blue edges. By the choice of $k$, $|A_p\cup R_G\cup X'|\le n \cdot 2^{k-1}$. But then $|B_G\cup X''|=|G|-|A_p\cup R_G\cup X'|\ge n \cdot 2^{k-1}+1$, a contradiction. Thus $G[R_G]$ has blue edges, as claimed. Then  $|R_G|\ge2$.  By Claim~\ref{R*_G},   $2\le |R_G|\le |R^*_G|\le 2n$.  

We first consider the case when  $|R^*_G| \ge n-1$. We claim that $|A_p\cup (X'\less R^*_G)|+|R^*_G|\le n \cdot 2^{k-2}+\max\{2n, k+n-1\}$. If $|R^*_G| \ge n$, then  by Lemma~\ref{n,n+1 Lemma}, $G[A_p]$ has no red edges and so $G[A_p\cup (X'\less R^*_G)]$ has no red edges.    By the choice of $k$, $|A_p\cup (X'\less R^*_G)|\le n \cdot 2^{k-2}$ and so $|A_p\cup (X'\less R^*_G)|+|R^*_G|\le n \cdot 2^{k-2}+2n$. If  $  |R_G^*| = n-1$, then applying Claim \ref{l-vertex}(ii) to   $R^*_G$ and $A_p$, $|A_p| \le n \cdot 2^{k-2}+2$. Thus 
$|A_p\cup (X'\less R^*_G)|+|R^*_G|\le n \cdot 2^{k-2}+2+(k-2) +(n-1)=n \cdot 2^{k-2} +k+n-1$.  Thus $|A_p\cup (X'\less R^*_G)|+|R^*_G|\le n \cdot 2^{k-2}+\max\{2n, k+n-1\}$, as claimed. But then
\[  
 |G|=|A_p\cup (X'\less R^*_G)|+|R^*_G|+|B_G\cup (X''\less R^*_G)| \le (n\cdot 2^{k-2}+\max\{2n, k+n-1\})+ n\cdot 2^{k-1}   < n \cdot 2^k+1,
\]
for all $k\ge3$, a contradiction. 
\medskip

It remains to consider the case   $ 2\le |R_G| \le  |R_G^*| \le n-2$. If   $|R_G^*| =3$, then $n=5$.  By applying Claim \ref{l-vertex}(iii) to   $R^*_G$ and $A_p$, $|A_p| \le 8 \cdot 2^{k-2}-1$.  But then 
 \[  
 |G|\le |A_p|+|B_G\cup X''|+|R^*_G|+|X'\less R^*_G|\le (8 \cdot 2^{k-2} - 1)+ 5\cdot 2^{k-1} +3+ (k-2) < 5 \cdot 2^k+1,
 \]
 for all $k \ge 3$, a contradiction. 
Thus  $|R^*_G|=|R_G|=2$. Then no vertex in $X$ is red-complete to $V(G)\less X$. Thus $|X''|\le k-2$.   Let $R_G=\{a, b\}$. Then $ab$ must be colored blue under $c$ because $G[R_G]$ has blue edges.  If $a$ or $b$, say $b$, is red-complete to $B_G$ in $G$, then neither $G[A_p\cup\{a\}\cup X']$ nor $G[B_G\cup\{b\}\cup X'']$ has blue edges. By minimality of $k$, $|A_p\cup \{a\}\cup X'|\le n \cdot 2^{k-1}$ and $|B_G\cup\{b\}\cup X''|\le n \cdot 2^{k-1}$. But then  $|G|=|A_p\cup \{a\}\cup X'|+|B_G\cup\{b\}\cup X''|\le n\cdot 2^{k-1}+ n\cdot 2^{k-1}< n\cdot 2^k+1$ for all $k\ge3$, a contradiction.  Thus neither $a$ nor $b$ is red-complete to $B_G$ in $G$. 
Let $a', b'\in B_G$ be such that $aa'$ and $bb'$ are colored blue under $c$. Then $a'=b'$, else we obtain a blue $C_{2n+1}$ in $G$ with vertices $a', a, b, b', x_1, y_1, x_2, \ldots,  y_{n-2},x_{n-1}$ in order, where $x_1, \ldots, x_{n-1}\in A_p$ and $y_1,\ldots, y_{n-2}\in B_G\less \{a', b'\}$, a contradiction.  Thus  $\{a, b\}$  is  red-complete to $B_G\less a'$ in $G$.  Then there exists $i\in[s]$ such that $A_i=\{a'\}$.   Since $G[B_G]$ has no blue edges, we see that 
 $\{a, b, a'\}$  must be  red-complete to $B_G\less a'$ in $G$.  By Claim~\ref{l-vertex}(ii,iii) applied to the three vertices $a, b, a'$ and  $B_G\less a'$,  we see that $|B_G\less a'|\le 4 \cdot 2^{k-2} +2$ when $n=4$ and $|B_G\less a'|\le 8 \cdot 2^{k-2} - 1$ when $n=5$. But then
\[
\begin{split}
 |G|	&= |A_p\cup X'|+|B_G\less a'|+|\{a, b, a'\}|+|X''|\\
  		& \le \begin{cases} 
 			4\cdot 2^{k-1} +(4\cdot 2^{k-2}+2)+3+(k-2),  &\text{when } n = 4\\
 			5\cdot 2^{k-1} +(8\cdot 2^{k-2}-1)+3+(k-2),  &\text{when } n=5\\
 		\end{cases}\\
 		& < n \cdot 2^k + 1
 \end{split}
 \]
 for all $k\ge3$, a contradiction.  Hence, $|A_{p-1}|\le n$.
 \qed\\
\end{pf}

By Claim~\ref{R*_G}, $ |R_G|\le  |R^*_G|\le 2n$.   We first consider the case when $ |R_G| \ge n$.   Since $|A_p|\ge n+1$, by Lemma~\ref{n,n+1 Lemma}, $G[A_p]$ has no red edges and no vertex in $X$ is red-complete to $V(G) \less X$. Thus  $|X| \le 2(k-1)$.  We first claim that $|B_G|\ge n$. Suppose $|B_G|\le n-1$.  
  If  $  |B_G| = n-1$,   then    $|A_p|\le n\cdot 2^{k-2}+2$ by Claim~\ref{l-vertex}(ii) applied to   $B_G$ and   $A_p$.  But then
  \[  
  |G|=|A_p|+|B_G|+|R_G|+|X|\le (n\cdot 2^{k-2}+2)+ (n-1)+2n+2(k-1)<n\cdot 2^k+1,
  \]
  for all  $k\ge3$,  a contradiction. Thus $3\le |B_G| \le n-2$. Then $n=5$ and  $  |B_G| =3$.  By Claim~\ref{l-vertex}(iii) applied to   $B_G$ and   $A_p$,   $|A_p|\le 8\cdot 2^{k-2}-1$.  But then 
  \[  
  |G|=|A_p|+|B_G|+|R_G|+|X|\le (8\cdot 2^{k-2}-1)+ 3+10+2(k-1)<5\cdot 2^k+1,
  \]
  for all  $k\ge3$,  a contradiction.  Thus $|B_G|\ge n$, as claimed.   
  By Lemma~\ref{n,n+1 Lemma}, $G[A_p]$ has no blue edges and no vertex in $X$ is blue-complete to $A_p$ in $G$.  Since $G[A_p\cup X']$ has neither red nor blue edges, and no vertex in $X$ is red- or blue-complete to $A_p$ in $G$, it follows that $|X''|\le k-2$  and  $|A_p\cup X'|\le n\cdot 2^{k-2}$ by minimality of $k$.   Then $|B_G|\ge n+1$, else   
  \[
  |G|=|A_p\cup X'|+|X''|+(|B_G|+|R_G|)\le n\cdot 2^{k-2}+(k-2)+( n+2n)<n\cdot 2^k+1,
  \]
  for all $k\ge3$, a contradiction. \medskip
  By  Lemma~\ref{n,n+1 Lemma}, $G[B_G]$ has no blue edges and so    $G[B_G\cup X'']$ has no blue edges. By minimality of $k$,   $|B_G\cup X''|\le n\cdot 2^{k-1}$.  But then 
  \[
  |G|=|A_p\cup X'|+|B_G\cup X''|+ |R_G|\le n\cdot 2^{k-2}+ n\cdot 2^{k-1}+2n<n\cdot 2^k+1,
  \]
  for all  $k\ge3$,  a contradiction.

  It remains to consider the case when $|R_G|\le n-1$.      Suppose first that $|B_G| \ge 2n+1$.  By Lemma \ref{n,n+1 Lemma}, $G[B_G]$ has no   blue edges. Thus all the edges in $\mathcal{R}[B]$ are colored red.  Since  $|A_{p-1}| \le n$  by Claim \ref{A_{p-1}}, we see that  $|B| \ge 3$, contrary to  Claim \ref{B}.
  Thus   $3\le |A_{p-1}|\le |B_G| \le 2n$.  If $|B_G| \ge  n-1$, by Claim \ref{l-vertex}(i,ii) applied to  $B_G$ and $A_p$  (and Lemma \ref{n,n+1 Lemma} applied to $B_G$ and $A_p$ to obtain $|X|\le 2(k-1)$ when $|B_G| \ge n$), we have
 \[
  |A_p| + |B_G|+|X|\le\begin{cases}  
(n \cdot 2^{k-1} +2) + (n-1)+2k, &  \text{if } |B_G| = n-1\\
n \cdot 2^{k-1} + 2n+2(k-1), &  \text{if } |B_G| \ge n.
 \end{cases}
 \]   
  Thus in either case, $|A_p|+|B_G|+|X|\le n \cdot 2^{k-1} + 2n+2k-2$. But then 
\[
|G|=(|A_p|+|B_G|+|X|)+|R_G|\le (n\cdot 2^{k-1}+ 2n+2k-2)+(n-1) <n\cdot 2^k+1, 
\]
for all $k \ge 3$, a contradiction.  Thus $3\le |B_G| \le n-2$. Then $ |B_G|=3$ and $n=5$.  If  $|R_G^*| \ge 4$ or $|B_G^*| \ge 4$, by applying Claim \ref{l-vertex}(ii) to any four vertices in $R_G^*$ or  $B_G^*$ and $A_p$,  we have  $|A_p| \le 5 \cdot 2^{k-1} +2$.  But  then
\[
|G|=|A_p|+|B_G|+|R_G|+|X|\le (5\cdot 2^{k-1}+ 2)+3+4+2k <5\cdot 2^k+1, 
\]
for all $k \ge 3$, a contradiction.
Thus  $|B_G|=|B_G^*| = 3$ and $|R_G| \le |R_G^*| \le 3$.  Then no vertex in $X$ is blue-complete to $V(G)\less X$. Thus $|X\less R_G^*|\le 2(k-2)$.  By Claim \ref{l-vertex}(iii) applied to   $B_G$ and $A_p$, $|A_p| \le 8 \cdot 2^{k-1} -1$. But then 
\[
  |G|=|A_p|+|B_G|+|R_G^*|+|X\less R_G^*|\le (8\cdot 2^{k-1}-1)+3+3+2(k-2) <5\cdot 2^k+1, 
\]
for all $k \ge 3$, a contradiction.  

This completes the proof of Theorem~\ref{C9C11}.\qed\\

  \section*{Acknowledgement}
The authors would like to thank Jingmei Zhang for  many helpful comments and discussion.

\end{document}